\documentclass[12pt,reqno]{amsart}

\usepackage{amsmath, epsfig, cite,color}
\usepackage{amssymb}
\usepackage{amssymb}
\usepackage{amsfonts}
\usepackage{latexsym}
\usepackage{graphicx}
\usepackage{algorithmic,algorithm}
\usepackage{amsthm,array,multirow,longtable}
\usepackage{tikz}

\newtheorem{thm}{Theorem}[section]
\newtheorem{prop}[thm]{Proposition}

\newtheorem{defi}[thm]{Definition}
\newtheorem{lem}[thm]{Lemma}
\newtheorem{en}[thm]{Entry}

\numberwithin{equation}{section}

\makeatletter \@addtoreset{equation}{section} \makeatother

\begin{document}
\allowdisplaybreaks
\title
{Four identities related to third order mock theta functions}
\author[S.-P. Cui]{Su-Ping Cui}
\address[S.-P. Cui]{Center for Combinatorics, LPMC, Nankai University, Tianjin 300071, P.R. China} \email{jiayoucui@163.com}
\author[N.S.S. Gu]{Nancy S.S. Gu}
\address[N.S.S. Gu]{Center for Combinatorics, LPMC, Nankai
University, Tianjin 300071, P.R. China} \email{gu@nankai.edu.cn}
\author[C.-Y. Su]{Chen-Yang Su}
\address[C.-Y. Su]{Center for Combinatorics, LPMC, Nankai
University, Tianjin 300071, P.R. China}
\email{cyangsu@163.com}

\date{\today}
\begin{abstract}
Ramanujan presented four identities for third order mock theta functions in his Lost Notebook.
In 2005, with the aid of complex analysis, Yesilyurt first proved these four identities. Recently, Andrews et al. provided different proofs by using $q$-series. In this paper, in view of some identities of a universal mock theta function
\begin{align*}
g(x;q)=x^{-1}\left(-1+\sum_{n=0}^{\infty}\frac{q^{n^{2}}}{(x;q)_{n+1}(qx^{-1};q)_{n}}\right),
\end{align*}
we establish new proofs of these four identities. In particular, by means of an identity of $g(x;q)$ given by Ramanujan and some theta function identities due to Mortenson, we find a new simple proof of the fourth identity.
\end{abstract}

\maketitle
\section{Introduction}

Let $q$ denote a complex number with $|q|<1$. Throughout this paper, we adopt the standard $q$-series notation \cite{Gasper-Rahman-2004}. For any positive integer $n$,
$$(a;q)_{0}:=1, \quad (a;q)_n:=\prod_{k=0}^{n-1}(1-aq^{k}), \quad (a;q)_{\infty}:=\prod_{k=0}^{\infty}(1-aq^{k}),$$
\begin{align}\label{jx}
j(x;q):=(x;q)_{\infty}(q/x;q)_{\infty}(q;q)_{\infty}.
\end{align}
Letting $a$ and $m$ be integers with $m$ positive, we define
\begin{align*}
J_{a,m}:=j(q^{a};q^{m}),\quad \quad \overline{J}_{a,m}:=j(-q^{a};q^{m}), \quad \quad  J_{m}:=J_{m,3m}=\prod_{k=1}^{\infty}(1-q^{mk}).
\end{align*}
Recall the Jacobi triple product identity \cite[Eq. (1.6.1)]{Gasper-Rahman-2004}
\begin{align*}
\sum_{k=-\infty}^{\infty}(-1)^{k}q^{\binom{k}{2}}x^{k}=(x;q)_{\infty}(q/x;q)_{\infty}(q;q)_{\infty}.
\end{align*}
Thus, Ramanujan's theta function $\varphi(q)$ and $\psi(q)$ are defined, respectively, as follows.
\begin{align}\label{varphiq}
\varphi(q)&:=\sum_{n=-\infty}^{\infty}q^{n^{2}}=(-q;q^{2})^{2}_{\infty}(q^{2};q^{2})_{\infty}
=\frac{J^{5}_{2}}{J^{2}_{1}J^{2}_{4}},\\
\psi(q)&:=\sum_{n=0}^{\infty}q^{n+1 \choose 2}=
\frac{(q^{2};q^{2})_{\infty}}{(q;q^{2})_{\infty}}=\frac{J^{2}_{2}}{J_{1}}.\label{4}
\end{align}
For any real number $a$, define
\begin{align}\label{3}
f_{a}(q):=\sum_{n=0}^{\infty}\frac{q^{n^{2}}}{(1+aq+q^{2})(1+aq^{2}+q^{4})\cdots (1+aq^{n}+q^{2n})}.
\end{align}
Set
\begin{align}\label{phi}
\tilde{\phi}(q):=f_0(q)=\sum_{n=0}^{\infty}\frac{q^{n^{2}}}{(-q^{2};q^{2})_{n}},
\end{align}
which is one of the third order mock theta functions included in Ramanujan's last letter to Hardy. Then Ramanujan \cite{ra} stated the following four identities related to third order mock theta functions in his Lost Notebook. For more on mock theta functions, see Andrews and Berndt's book\cite{part5}.

\begin{en}\cite[p. 2]{ra}\label{401}
Suppose that $a$ and $b$ are real numbers such that $a^{2}+b^{2}=4$. Recall that $f_{a}(q)$ is
defined by (\ref{3}). Then
\begin{align}\label{41} \nonumber
&\frac{b-a+2}{4}f_{a}(-q)+\frac{b+a+2}{4}f_{-a}(-q)-\frac{b}{2}f_{b}(q)\\
&\quad=\frac{(q^{4};q^{4})_{\infty}}{(-q;q^{2})_{\infty}}\prod_{n=1}^{\infty}\frac{1-bq^{n}+q^{2n}}
{1+(a^{2}b^{2}-2)q^{4n}+q^{8n}}.
\end{align}
\end{en}

\begin{en}\cite[p. 2]{ra}\label{402}
Let $a$ and $b$ be real numbers with $a^{2}+ab+b^{2}=3$. Then, with $f_{a}(q)$
defined by (\ref{3}),
\begin{align}\label{42} \nonumber
&(a+1)f_{-a}(q)+(b+1)f_{-b}(q)-(a+b-1)f_{a+b}(q)\\
&\quad=3\frac{(q^{3};q^{3})_{\infty}^{2}}{(q;q)_{\infty}}\prod_{n=1}^{\infty}\frac{1}
{1+ab(a+b)q^{3n}+q^{6n}}.
\end{align}
\end{en}

\begin{en}\cite[p. 17]{ra}\label{403}
Let $\psi(q)$ and $f_{a}(q)$ be defined by \eqref{4} and \eqref{3}, respectively. Then
\begin{align*}
&\frac{1+\sqrt{3}}{2}f_{-1}(-q)+\frac{3+\sqrt{3}}{6}f_{1}(-q)-f_{\sqrt{3}}(q)\\
&\quad=\frac{2}{\sqrt{3}}\psi(-q)\frac{(q^{4};q^{4})_{\infty}}{(q^{6};q^{6})_{\infty}}
\prod_{n=1}^{\infty}\frac{1}{1+\sqrt{3}q^{n}+q^{2n}}.
\end{align*}
\end{en}

\begin{en}\cite[p. 17]{ra}\label{404}
Let $\psi(q)$, $f_{a}(q)$, and $\tilde{\phi}(q)$ be defined by \eqref{4}, \eqref{3}, and \eqref{phi}, respectively. Then
\begin{align*}
&\frac{1}{2}(1+e^{\pi i/4})\tilde{\phi}(iq)+\frac{1}{2}(1+e^{-\pi i/4})\tilde{\phi}(-iq)-f_{\sqrt{2}}(q)\\
&\quad=\frac{1}{\sqrt{2}}\psi(-q)(-q^{2};q^{4})_{\infty}
\prod_{n=1}^{\infty}\frac{1}{1+\sqrt{2}q^{n}+q^{2n}}.
\end{align*}
\end{en}

In 1944, Dyson \cite{dyson} defined the rank of a partition to be the largest part minus the number of parts. Let $N(m,n)$ denote the number of partitions of the positive integer $n$ with rank $m$. Then he showed that the generating function for $N(m,n)$ is given by
\begin{align}\label{2}
\sum_{m=-\infty}^{\infty}\sum_{n=0}^{\infty}N(m,n)q^{n}x^{m}=\sum_{n=0}^{\infty}
\frac{q^{n^{2}}}{(qx;q)_{n}(qx^{-1};q)_{n}}=:G(x,q).
\end{align}
From \eqref{3} and \eqref{2}, it can be seen that
\begin{align*}
f_{a}(q)=G\left(\frac{-a\pm \sqrt{a^{2}-4}}{2},q\right).
\end{align*}

The proofs of Entries 1.1-1.4 which were first provided by Yesilyurt \cite{ye} rely on the following lemma due to Atkin and Swinnerton-Dyer.
\begin{lem}\cite{atkindyer} Let $q$, $|q|< 1$, be fixed. Suppose that $\vartheta(z)$ is an analytic
function of $z$, except for possibly a finite number of poles, in every region,
$0 < z_{1}\leq |z| \leq z_{2}$. If
$$\vartheta(zq) = Az^k\vartheta(z)$$
for some integer $k$ (positive, zero, or negative) and some constant $A$, then
either $\vartheta(z)$ has $k$ more poles than zeros in the region $|q|< |z|\leq 1$,
or $\vartheta(z)$ vanishes identically.
\end{lem}
Recently, Andrews et al. \cite{an} offered different proofs of these four identities by using $q$-series. The proofs of Entries 1.1-1.3 depend on a partial fraction expansion formula, and in the proof of Entry 1.4, two $2$-dissentions for two special cases of the rank generating function $G(x,q)$ are the main tools.

The object of this paper is to present new proofs of these four identities. Recall a universal mock theta function,
\begin{align}\label{20}
g(x;q):=x^{-1}\left(-1+\sum_{n=0}^{\infty}\frac{q^{n^{2}}}{(x;q)_{n+1}(qx^{-1};q)_{n}}\right).
\end{align}
A transformation formula of $g(x;q)$ given by Ramanujan \cite{ra} and an identity of $g(x;q)$ due to Hickerson \cite{hickson,hm2} play important roles in our proofs, see Propositions \ref{25} and \ref{gm} in the next section. In particular, in view of some theta function identities due to Mortenson \cite{mo,mo2}, we establish a simple proof of Entry 1.4. This paper is organized as follows. In Section 2, we show some results which are used in our proofs.
In Section 3, we prove Entries 1.1-1.4.

\section{Preliminaries}

In this section, we present some preliminary results. For later use, we need the following identities:
\begin{align*}
\overline{J}_{1,2}
=\frac{J^5_2}{J^2_1J^2_4}, \quad
J_{1,2}=\frac{J^2_1}{J_2},
\quad J_{1,4}=\frac{J_1J_4}{J_2}.
\end{align*}
In addition, the following general identities are frequently used in this paper.
\begin{align}\label{xq}
j(qx;q)=-x^{-1}j(x;q),
\end{align}
\begin{align}\label{q/x}
j(x;q)=j(qx^{-1};q),
\end{align}
\begin{align}
j(x^{2};q^{2})=j(x;q)j(-x;q)/J_{1,2},\label{x2q2}
\end{align}
\begin{align}
j(x;q)=J_{1}j(x;q^{2})j(qx;q^{2})/J^{2}_{2}, \label{xqx}
\end{align}
\begin{align}
j(x;q)=\sum^{m-1}_{k=0}(-1)^kq^{{k\choose 2}}x^kj\left((-1)^{m+1}q^{{m\choose2}+mk}x^m;q^{m^2}\right)\label{mz}.
\end{align}
Setting $m=2$ in \eqref{mz} yields that
\begin{align} \label{22}
j(x;q)=j(-qx^{2};q^{4})-xj(-q^{3}x^{2};q^{4}).
\end{align}

Hickerson and Mortenson \cite{hm2} defined Appell-Lerch sums as follows.
\begin{defi}
Let $x,z\in\mathbb{C}^*: =\mathbb{C}\backslash \{0\}$ with neither $z$ nor $xz$ an integral power of $q$. Then
\begin{align*}
m(x,q,z):=\frac{1}{j(z;q)}\sum_{r=-\infty}^{\infty}\frac{(-1)^rq^{\binom{r}{2}}z^r}{1-q^{r-1}xz}.
\end{align*}
\end{defi}
Following \cite{hm2}, the term ``generic" means that the parameters do not cause poles in the Appell-Lerch sums or in the quotients of theta functions. The next proposition can be found in the Lost Notebook.
\begin{prop}\label{25}(\cite[p. 32]{ra}, \cite[Eq. (12.5.13)]{part1}). For generic $x\in \mathbb{C}^*$,
\begin{align*}
g(x;q)=-x^{-1}+qx^{-3}g(-qx^{-2};q^{4})-qg(-qx^{2};q^{4})+\frac{J_{2}J_{2,4}^{2}}{xj(x;q)j(-qx^{2};q^{2})}.
\end{align*}
\end{prop}
The following proposition was first proved by Hickerson \cite{hickson}, and then Hickerson and Mortenson \cite{hm2} rewrote the identity in terms of Appell-Lerch sums.
\begin{prop}\label{gm}(\cite[Theorem 2.2]{hickson}, \cite{hm2}). For generic $x, z\in \mathbb{C}^*$,
\begin{align}
g(x;q)=-x^{-2}m(qx^{-3},q^{3},x^{3}z)-x^{-1}m(q^{2}x^{-3},q^{3},x^{3}z)+
\frac{J_{1}^{2}j(xz;q)j(z;q^{3})}{j(x;q)j(z;q)j(x^{3}z;q^{3})}. \label{gms}
\end{align}
\end{prop}
Notice that a special case of \eqref{gms} with $z=1$ reads
\begin{align}
g(x;q)=-x^{-2}m(qx^{-3},q^{3},x^{3})-x^{-1}m(q^{2}x^{-3},q^{3},x^{3})+
\frac{J_{3}^{3}}{J_{1}j(x^{3};q^{3})}.\label{z=1}
\end{align}

In the proof of Entry 1.4, we apply the next two propositions to simplify theta function identities.
\begin{prop}\label{23}\cite[Proposition 3.4]{mo2} Let $x\neq 0$. Then
\begin{align*}
j(q^{2}x;q^{4})j(q^{5}x;q^{8})+\frac{q}{x}\cdot j(x;q^{4})j(qx;q^{8})-
\frac{J_{1}}{J_{4}}\cdot j(-q^{3}x;q^{4})j(q^{3}x;q^{8})=0.
\end{align*}
\end{prop}
\begin{prop}\label{24}\cite[Proposition 2.5]{mo} Let $x\neq 0$. Then
\begin{align*}
j(-x;q^{4})j(-q^{5}x;q^{8})-j(-q^{2}x;q^{4})j(-qx;q^{8})-
x\frac{J_{1}}{J_{4}}\cdot j(q^{3}x;q^{4})j(-q^{7}x;q^{8})=0.
\end{align*}
\end{prop}

Furthermore, in view of \eqref{2} and \eqref{20}, we arrive at
\begin{align}\label{gg}
G(x,q)=(1-x)\left(xg(x;q)+1\right).
\end{align}

\section{Proofs of Entries 1.1-1.4}
In this section, employing Propositions \ref{25}-\ref{24}, we prove Entries 1.1-1.4.

\noindent{\textbf{\emph{Proof of Entry \ref{401}.}}}
Set $a=2\cos\theta$, $b=2\sin\theta$, and $t=e^{i\theta}$. Then $a=t+t^{-1}$ and $b=-i(t-t^{-1})$.
So we have
\begin{align*}
&\frac{b-a+2}{4}=\frac{(i-1)(1-it)(1-t)}{4t},\\
&\frac{b+a+2}{4}=\frac{(i+1)(1-it)(1+t)}{4t},\\
&\frac{b}{2}=\frac{i(1-t^{2})}{2t}, \quad a^{2}b^{2}-2=-(t^{4}+t^{-4}).
\end{align*}
Moreover, according to \eqref{3} and \eqref{2}, we obtain that
$$f_{a}(q)=G(-t,q),\quad f_{-a}(q)=G(t,q),\quad f_{b}(q)=G(it,q).$$ Thus, using the above relations, we  rewrite \eqref{41} as
\begin{align}\label{410}\nonumber
&\frac{(i-1)(1-it)(1-t)}{4t}G(-t,-q)+\frac{(i+1)(1-it)(1+t)}{4t}G(t,-q)\\
&\qquad -\frac{i(1-t^{2})}{2t}G(it,q)
=\frac{(q^{4};q^{4})_{\infty}(-iqt;q)_{\infty}(iqt^{-1};q)_{\infty}}{(-q;q^{2})_{\infty}
(q^{4}t^{4};q^{4})_{\infty}(q^{4}t^{-4};q^{4})_{\infty}}.
\end{align}
Multiplying both sides of \eqref{410} by $4t(1+it)/((i-1)(1-t^4))$, we arrive at
\begin{align}\label{411}
\frac{1}{1+t}G(-t,-q)-\frac{i}{1-t}
G(t,-q)+\frac{i-1}{1-it}G(it,q)=-\frac{2(1+i)tJ^{3}_{4}j(-it;q)}{J^{2}_{2}j(t^4;q^4)}.
\end{align}
Hence, to prove Entry 1.1, it suffices to prove \eqref{411}.

First, we consider the left-hand side of \eqref{411}. Using \eqref{gg} in the first equality and
Proposition \ref{25} in the third equality below, we have
\begin{align}
&\frac{1}{1+t}G(-t,-q)-\frac{i}{1-t}G(t,-q)+\frac{i-1}{1-it}G(it,q)\nonumber \\
\quad&=-tg(-t;-q)+1-i(tg(t;-q)+1)+(i-1)(itg(it;q)+1) \nonumber\\
&=-tg(-t;-q)-itg(t;-q)-(1+i)tg(it;q) \nonumber\\
&=-1-qt^{-2}g(qt^{-2};q^4)
-qtg(qt^2;q^4)+\frac{J_2J^2_{2,4}}{j(-t;-q)j(qt^2;q^2)}\nonumber \\
&\quad+i+iqt^{-2}g(qt^{-2};q^4)
-iqtg(qt^2;q^4)-\frac{i J_2J^2_{2,4}}{j(t;-q)j(qt^2;q^2)}\nonumber \\
&\quad-i+1+(1-i)qt^{-2}g(qt^{-2};q^4)
+(1+i)qtg(qt^2;q^4)-\frac{(1-i)J_2J^2_{2,4}}{j(it;q)j(qt^2;q^2)}\nonumber \\
&=(-1+i+1-i)qt^{-2}g(qt^{-2};q^4)
-(1+i-1-i)qtg(qt^2;q^4)\nonumber \\
&\quad+\frac{J_{2}J_{2,4}^{2}}{j(-t;-q)j(qt^{2};q^{2})}-\frac{i J_{2}J_{2,4}^{2}}{j(t;-q)j(qt^{2};q^{2})}
-\frac{(1-i)J_{2}J_{2,4}^{2}}{j(it;q)j(qt^{2};q^{2})}\nonumber \\
&=\frac{J_{2}J_{2,4}^{2}}{j(qt^{2};q^{2})}\left(\frac{1}{j(-t;-q)}-\frac{i}{j(t;-q)}-
\frac{1-i}{j(it;q)}\right).\label{lhs}
\end{align}
Define
\begin{align*}
A&:=\frac{1}{j(-t;-q)}-\frac{i}{j(t;-q)}-\frac{1-i}{j(it;q)}\\
&=\frac{j(t;-q)j(it;q)-ij(-t;-q)j(it;q)-(1-i)j(-t;-q)j(t;-q)}
{j(-t;-q)j(t;-q)j(it;q)}.
\end{align*}
To evaluate $A$, by \eqref{22}, we deduce that
\begin{align*}
j(-t;-q)&=j(qt^{2};q^{4})+tj(q^{3}t^{2};q^{4}),\\
j(t;-q)&=j(qt^{2};q^{4})-tj(q^{3}t^{2};q^{4}),\\
j(it;q)&=j(qt^{2};q^{4})-itj(q^{3}t^{2};q^{4}).
\end{align*}
Then by means of the above three identities, we find that
\begin{align}
&j(t;-q)j(it;q)-ij(-t;-q)j(it;q)-(1-i)j(-t;-q)j(t;-q) \nonumber \\
&=\left(j(qt^{2};q^{4})-tj(q^{3}t^{2};q^{4})\right)\left(j(qt^{2};q^{4})-itj(q^{3}t^{2};q^{4})\right) \nonumber \\
&\quad-i\left(j(qt^{2};q^{4})+tj(q^{3}t^{2};q^{4})\right)\left(j(qt^{2};q^{4})-itj(q^{3}t^{2};q^{4})\right) \nonumber \\
&\quad-(1-i)\left(j(qt^{2};q^{4})+tj(q^{3}t^{2};q^{4})\right)\left(j(qt^{2};q^{4})-tj(q^{3}t^{2};q^{4})\right) \nonumber \\
&=-2(1+i)tj(qt^{2};q^{4})j(q^{3}t^{2};q^{4}).\label{A}
\end{align}
Hence, utilizing \eqref{x2q2}, \eqref{xqx}, and \eqref{A}, we obtain that
\begin{align*}
A&=-\frac{2(1+i)tj(qt^{2};q^{4})j(q^{3}t^{2};q^{4})}{j(-t;-q)j(t;-q)j(it;q)}\\
&=-\frac{2(1+i)tJ_{4}^{2}j(qt^{2};q^{2})}{J_{2}j(-t;-q)j(t;-q)j(it;q)}\cdot\frac{j(-it;q)}{j(-it;q)}\\
&=-\frac{2(1+i)tJ_{4}^{2}j(qt^{2};q^{2})j(-it;q)}{J_{2}J_{1,2}\overline{J}_{1,2}J_{2,4}j(t^4;q^4)}.
\end{align*}
Thus, the right-hand side of \eqref{lhs} becomes
\begin{align*}
&\frac{J_{2}J_{2,4}^{2}}{j(qt^{2};q^{2})}\left(-\frac{2(1+i)tJ_{4}^{2}j(qt^{2};q^{2})j(-it;q)}
{J_{2}J_{1,2}\overline{J}_{1,2}J_{2,4}j(t^4;q^4)}\right) \\
&=-\frac{2(1+i)tJ^{3}_{4}j(-it;q)}{J^{2}_{2}j(t^4;q^4)},
\end{align*}
which is the right-hand side of \eqref{411}. Therefore, we complete the proof of Entry 1.1. \qed

\noindent{\textbf{\emph{Proof of Entry \ref{402}.}}} We parameterize $a^2+ab+b^2=3$ by $a=2\cos(\theta+2\pi/3)$, $b=2\cos\theta$, and $t=e^{i\theta}$. Let $\omega=e^{2\pi i/3}$. Then we have $a=\omega t+(\omega t)^{-1}$ and $b=t+t^{-1}$. So,
$$a+b=-\omega^{2}t-(w^{2}t)^{-1},\quad a+1=\frac{1-t^{3}}{\omega t(1-\omega t)},\quad b+1=\frac{1-t^{3}}{t(1-t)},$$
$$a+b-1=-\frac{1-t^{3}}{\omega^{2}t(1-\omega^{2}t)}, \quad ab(a+b)=-t^{3}-t^{-3},$$
and
$$f_{-a}(q)=G(\omega t,q),\quad f_{-b}(q)=G(t,q), \quad f_{a+b}(q)=G(\omega^{2} t,q).$$
Therefore, \eqref{42} is equivalent to the following identity:
\begin{align}
&\frac{1-t^{3}}{\omega t(1-\omega t)}G(\omega t,q)+\frac{1-t^{3}}{t(1-t)}G(t,q)+
\frac{1-t^{3}}{\omega^{2}t(1-\omega^{2}t)}G(\omega^{2}t,q)\nonumber \\
&\qquad=\frac{3(q^{3};q^{3})_{\infty}^{2}}{(q;q)_{\infty}}\prod_{n=1}^{\infty}
\frac{1}{1-(t^{3}+t^{-3})q^{3n}+q^{6n}}.\label{G123}
\end{align}
Observing that the right-hand side of \eqref{G123} reduces to
$$\frac{3J^{3}_{3}(1-t^3)}{J_{1}j(t^3;q^3)},$$
we are required to prove that
\begin{align*}
&\frac{G(\omega t,q)}{\omega t(1-\omega t)}+\frac{G(t,q)}{t(1-t)}+
\frac{G(\omega^{2}t,q)}{\omega^{2}t(1-\omega^{2}t)}
=\frac{3J^3_3}{J_1j(t^3;q^3)}.
\end{align*}
From \eqref{gg}, it can be seen that
\begin{align}\label{421} \nonumber
&\frac{G(\omega t,q)}{\omega t(1-\omega t)}+\frac{G(t,q)}{t(1-t)}+
\frac{G(\omega^{2}t,q)}{\omega^{2}t(1-\omega^{2}t)}\\ \nonumber
&=g(\omega t;q)+\frac{1}{\omega t}+g(t;q)+\frac{1}{t}
+g(\omega^{2}t;q)+\frac{1}{\omega^{2}t}\\
&=g(t;q)+g(\omega t;q)+g(\omega^{2}t;q).
\end{align}
Then applying \eqref{z=1} to the right-hand side of \eqref{421} gives
\begin{align*}
&\frac{G(\omega t,q)}{\omega t(1-\omega t)}+\frac{G(t,q)}{t(1-t)}+
\frac{G(\omega^{2}t,q)}{\omega^{2}t(1-\omega^{2}t)}\\
=&-t^{-2}m(qt^{-3},q^{3},t^{3})-t^{-1}m(q^{2}t^{-3},q^{3},t^{3})+
\frac{J_{3}^{3}}{J_{1}j(t^{3};q^{3})}
\\&-\omega t^{-2}m(qt^{-3},q^{3},t^{3})-\omega^{2}t^{-1}m(q^{2}t^{-3},q^{3},t^{3})+
\frac{J_{3}^{3}}{J_{1}j(t^{3};q^{3})}\\
&-\omega^{2}t^{-2}m(qt^{-3},q^{3},t^{3})-\omega t^{-1}m(q^{2}t^{-3},q^{3},t^{3})+
\frac{J_{3}^{3}}{J_{1}j(t^{3};q^{3})}\\
=&-(1+\omega+\omega^{2})t^{-2}m(qt^{-3},q^{3},t^{3})-(1+\omega+\omega^{2})t^{-1}m(q^{2}t^{-3},q^{3},t^{3})\\
&+\frac{3J^3_3}{J_1j(t^3;q^3)}\\
=&\frac{3J^3_3}{J_1j(t^3;q^3)}.
\end{align*}
Therefore, we complete the proof of Entry 1.2.
\qed

\noindent{\textbf{\emph{Proof of Entry \ref{403}.}}}
It is easy to find that Entry 1.3 is a special case of Entry 1.1. Letting $a=1$ and $b=\sqrt{3}$ in \eqref{41}, we derive that
\begin{align*}
\frac{1+\sqrt{3}}{4}f_{1}(-q)+\frac{3+\sqrt{3}}{4}f_{-1}(-q)-\frac{\sqrt{3}}{2}f_{\sqrt{3}}(q)
=\frac{J_{1}J^{2}_{4}}{J^{2}_{2}}\prod_{n=1}^{\infty}\frac{1-\sqrt{3}q^{n}+q^{2n}}{1+q^{4n}+q^{8n}},
\end{align*}
which is equivalent to the following identity:
\begin{align}\label{432}
\frac{3+\sqrt{3}}{6}f_{1}(-q)+\frac{1+\sqrt{3}}{2}f_{-1}(-q)-f_{\sqrt{3}}(q)
=\frac{2J_{1}J^{2}_{4}}{\sqrt{3}J^{2}_{2}}\prod_{n=1}^{\infty}\frac{1-\sqrt{3}q^{n}+q^{2n}}
{1+q^{4n}+q^{8n}}.
\end{align}
Then the right-hand side of \eqref{432} becomes
\begin{align*}
&\frac{2J_{1}J^{2}_{4}}{\sqrt{3}J^{2}_{2}}\prod_{n=1}^{\infty}\frac{(1-\sqrt{3}q^{n}+q^{2n})(1+\sqrt{3}q^{n}+q^{2n})}
{(1+q^{4n}+q^{8n})(1+\sqrt{3}q^{n}+q^{2n})}\\
&=\frac{2J_{1}J^{2}_{4}}{\sqrt{3}J^{2}_{2}}\prod_{n=1}^{\infty}\frac{1-q^{2n}+q^{4n}}
{1+q^{4n}+q^{8n}}\prod_{n=1}^{\infty}\frac{1}{1+\sqrt{3}q^{n}+q^{2n}}\\
&=\frac{2J_{1}J^{2}_{4}}{\sqrt{3}J^{2}_{2}}\prod_{n=1}^{\infty}\frac{(1+q^{6n})(1-q^{4n})}
{(1+q^{2n})(1-q^{12n})}\prod_{n=1}^{\infty}\frac{1}{1+\sqrt{3}q^{n}+q^{2n}}\\
&=\frac{2J_{1}J^{2}_{4}}{\sqrt{3}J^{2}_{2}}\cdot \frac{J_2}{J_6}\prod_{n=1}^{\infty}\frac{1}{1+\sqrt{3}q^{n}+q^{2n}}\\
&=\frac{2}{\sqrt{3}}\psi(-q)\frac{J_4}{J_6}\prod_{n=1}^{\infty}\frac{1}{1+\sqrt{3}q^{n}+q^{2n}},
\end{align*}
where we use \eqref{4} in the last equality. Therefore, combining \eqref{432} and the above identity, we complete the proof of Entry 1.3.
\qed

\noindent{\textbf{\emph{Proof of Entry \ref{404}.}}}
Let $\alpha=e^{\pi i/4}$. According to \eqref{3}, \eqref{phi}, and \eqref{2}, we have
\begin{align*}
f_{\sqrt{2}}(q)=G(-\alpha,q),\quad \tilde{\phi}(q)=G(i,q).
\end{align*}
In view of the above two identities, we restate Entry \ref{404} as
\begin{align*}
\frac{1+\alpha}{2}G(i,iq)+\frac{1+\alpha^{-1}}{2}G(i,-iq)-G(-\alpha,q)=
\frac{1}{\sqrt{2}}\frac{\psi(-q)(-q^{2};q^{4})_{\infty}}{(-q\alpha ;q)_{\infty}(-q\alpha^{-1};q)_{\infty}}.
\end{align*}
Dividing both sides of the above identity by $(1+\alpha)/2$, we deduce that
\begin{align}\label{441}
G(i,iq)+\alpha^{-1}G(i,-iq)-\frac{2}{1+\alpha}G(-\alpha,q)=
\sqrt{2}\frac{\psi(-q)(-q^{2};q^{4})_{\infty}(q;q)_{\infty}}{j(-\alpha;q)}.
\end{align}
Replacing $q$ by $iq$ in \eqref{441} yields that
\begin{align}\label{39}
G(i,-q)+\alpha^{-1}G(i,q)-\frac{2}{1+\alpha}G(-\alpha,iq)
=\sqrt{2}\frac{\psi(-iq)(q^{2};q^{4})_{\infty}(iq;iq)_{\infty}}{j(-\alpha;iq)}.
\end{align}
Therefore, to obtain Entry 1.4, it suffices to prove \eqref{39}.

In the following, we first consider the left-hand side of \eqref{39}. Using \eqref{gg} and
Proposition \ref{25}, we have
\begin{align}\label{442} \nonumber
G(i,q)&=(i+1)g(i;q)+1-i\\
&=(i+1)\left(i+(i-1)qg(q;q^{4})+\frac{J_{2}J_{2,4}^{2}}{ij(i;q)j(q;q^{2})}\right)+1-i\nonumber\\
&=-2qg(q;q^{4})+(1-i)\frac{J_{2}J_{2,4}^{2}}{j(i;q)j(q;q^{2})}.
\end{align}
Similarly, by noticing that $\alpha^2=i$ and $\alpha^3=-\alpha^{-1}$, we arrive at
\begin{align}
G(i,-q)&=2qg(-q;q^{4})+(1-i)\frac{J_{2}J_{2,4}^{2}}{j(i;-q)j(-q;q^{2})}
\label{443},\\
-\frac{2}{1+\alpha}G(-\alpha,q)&=2iqg(iq;q^{4})-2q\alpha g(-iq;q^{4})-\frac{2J_{2}J_{2,4}^{2}}{j(-\alpha;q)j(-iq;q^{2})}\label{444}.
\end{align}
To evaluate the right-hand side of \eqref{444}, using \eqref{jx} and \eqref{x2q2}, we find that
\begin{align}
j(\alpha;q)j(-\alpha;q)
&=J_{1,2}j(i;q^2)=(1-i)\frac{J^{2}_{1}J_{8}}{J_{4}},\label{ja}
\end{align}
\begin{align}\label{ji}
j(-iq;q^{2})=\frac{J^{2}_{4}}{J_{8}}.
\end{align}
So, by \eqref{varphiq}, \eqref{ja}, and \eqref{ji}, we deduce that
\begin{align*}
\frac{J_2J^2_{2,4}}{j(-\alpha;q)j(-iq;q^2)}
=\frac{J_2J^2_{2,4}}{j(-\alpha;q)j(-iq;q^2)}\cdot\frac{j(\alpha;q)}{j(\alpha;q)}
=\frac{J_{2}J^2_{2,4}j(\alpha;q)}{(1-i)J^2_1J_{4}}=\frac{\varphi(q)j(\alpha;q)}{(1-i)J_4}.
\end{align*}
Therefore, we restate \eqref{444} as
\begin{align*}
-\frac{2}{1+\alpha}G(-\alpha,q)&=2iqg(iq;q^{4})-2q\alpha g(-iq;q^{4})-\frac{2\varphi(q)j(\alpha;q)}{(1-i)J_4}.
\end{align*}
Setting $q$ by $iq$ in the above identity yields that
\begin{align}
-\frac{2}{1+\alpha}G(-\alpha,iq)&=-2qg(-q;q^{4})+2q\alpha^{-1} g(q;q^{4})-\frac{2\varphi(iq)j(\alpha;iq)}{(1-i)J_{4}}.\label{qiq}
\end{align}
Then substituting \eqref{442}, \eqref{443}, and \eqref{qiq} into the left-hand side of \eqref{39}, we obtain
\begin{align}\label{445} \nonumber
&G(i,-q)+\alpha^{-1}G(i,q)-\frac{2}{1+\alpha}G(-\alpha,iq)\\ \nonumber
&=(1-i)\frac{J_{2}J_{2,4}^{2}}{j(i;-q)j(-q;q^{2})}+(1-i)\alpha^{-1}\frac{J_{2}J_{2,4}^{2}}{j(i;q)j(q;q^{2})}
-\frac{2\varphi(iq)j(\alpha;iq)}{(1-i)J_{4}}\\
&=(1-i)\frac{J_{2}J_{2,4}^{2}\left(j(i;q)j(q,q^{2})+\alpha^{-1}j(i;-q)j(-q;q^{2})\right)}
{j(i;-q)j(-q;q^{2})j(i;q)j(q;q^{2})}-\frac{2\varphi(iq)j(\alpha;iq)}{(1-i)J_{4}}.
\end{align}
To continue the evaluation from \eqref{445}, we note that from \eqref{q/x} and \eqref{22},
$$j(i;q)=j(q;q^4)-ij(q^3;q^4)=(1-i)j(q;q^{4}).$$
Substituting the above identity into \eqref{445}, we have
\begin{align}\label{446}
&G(i,-q)+\alpha^{-1}G(i,q)-\frac{2}{1+\alpha}G(-\alpha,iq)\nonumber \\
&=\frac{J_{2}J_{2,4}^{2}(j(q;q^4)j(q;q^{2})+\alpha^{-1}j(-q;q^4)j(-q;q^{2}))}
{j(-q;q^4)j(q;q^4)j(-q;q^2)j(q;q^2)}-\frac{2\varphi(iq)j(\alpha;iq)}{(1-i)J_{4}}\nonumber \\
&=\frac{j(q;q^{4})j(q;q^{2})+\alpha^{-1}j(-q;q^{4})j(-q;q^{2})}{J_{4}}-\frac{2\varphi(iq)j(\alpha;iq)}{(1-i)J_{4}},
\end{align}
where we use \eqref{x2q2} in the last equality.

Next, we give an alternative representation for the right-hand side of \eqref{39}.
\begin{align}\label{rhs3.9} \nonumber
&\frac{\sqrt{2}\psi(-iq)(q^{2};q^{4})_{\infty}(iq;iq)_{\infty}}{j(-\alpha;iq)}
\nonumber\\
&=\frac{\sqrt{2}\psi(-iq)(q^{2};q^{4})_{\infty}(iq;iq)_{\infty}}{j(-\alpha;iq)}
\cdot \frac{j(\alpha;iq)}{j(\alpha;iq)}\nonumber\\
&=\frac{\sqrt{2}J_{2,4}j(\alpha;iq)}{(1-i)J_{4}},
\end{align}
where we utilize \eqref{4} and \eqref{ja}. Therefore, combining
\eqref{446} and \eqref{rhs3.9}, we only need to prove that
\begin{align}\label{449}
(1-i)\left(j(q;q^{4})j(q;q^{2})+\alpha^{-1}j(-q;q^{4})j(-q;q^{2})\right)-
2\varphi(iq)j(\alpha;iq)=\sqrt{2}j(q^2;q^4)j(\alpha;iq).
\end{align}
In view of \eqref{22}, we have
\begin{align}
&j(q;q^{4})=j(-q^{6};q^{16})-qj(-q^{14};q^{16}), \label{qq4}\\
&j(q;q^{2})=j(-q^{4};q^{8})-qj(-q^{8};q^{8}),\\
&\varphi(iq)=j(-iq;-q^2)=j(-q^{4};q^{8})+iqj(-q^{8};q^{8}).
\end{align}
Moreover, using \eqref{22} twice yields
\begin{align}\label{jq}
j(\alpha;iq)&=j(q;q^{4})-\alpha j(-q;q^{4})\nonumber\\
&=j(-q^{6};q^{16})-qj(-q^{14};q^{16})
-\alpha\left(j(-q^{6};q^{16})+qj(-q^{14};q^{16})\right)\nonumber\\
&=(1-\alpha)j(-q^{6};q^{16})-(1+\alpha)q j(-q^{14};q^{16}).
\end{align}
Then substituting \eqref{qq4}-\eqref{jq} into the left-hand side of \eqref{449}, we find that
\begin{align}\label{4410} \nonumber
&(1-i)\left(j(q;q^{4})j(q;q^{2})+\alpha^{-1}j(-q;q^{4})j(-q;q^{2})\right)-
2\varphi(iq)j(\alpha;iq)\nonumber\\
&=(1-i)\left(j(-q^6;q^{16})-qj(-q^{14};q^{16})\right)
\left(j(-q^4;q^{8})-qj(-q^{8};q^{8})\right)\nonumber\\
&\quad+(1-i)\alpha^{-1}\left(j(-q^6;q^{16})+qj(-q^{14};q^{16})\right)
\left(j(-q^4;q^{8})+qj(-q^{8};q^{8})\right)\nonumber\\
&\quad-2\left(j(-q^4;q^{8})+iqj(-q^{8};q^{8})\right)
\left((1-\alpha)j(-q^6;q^{16})-(1+\alpha)qj(-q^{14};q^{16})\right)\nonumber\\
&=(-1+\sqrt{2}-i)\left(j(-q^{4};q^{8})j(-q^{6};q^{16})-q^{2}j(-q^{8};q^{8})j(-q^{14};q^{16})\right)\nonumber\\
&\quad+(1+\sqrt{2}+i)q\left(j(-q^{4};q^{8})j(-q^{14};q^{16})-j(-q^{8};q^{8})j(-q^{6};q^{16})\right).
\end{align}
Furthermore, setting $q\rightarrow q^{2}$ and $x=-1$ in Proposition \ref{23} and employing \eqref{q/x}, we arrive at
\begin{align}\label{4411}
&j(-q^{4};q^{8})j(-q^{6};q^{16})-q^{2}j(-q^{8};q^{8})j(-q^{14};q^{16})
=j(q^2;q^4)j(-q^6;q^{16}).
\end{align}
Setting $q\rightarrow q^{2}$ and $x=q^{4}$ in Proposition \ref{24} and applying \eqref{xq}
and \eqref{q/x}, we obtain
\begin{align}
j(-q^{4};q^{8})j(-q^{14};q^{16})-j(-q^{8};q^{8})j(-q^{6};q^{16})=-j(q^2;q^4)j(-q^{14};q^{16}).\label{4412}
\end{align}
Then substituting \eqref{4411} and \eqref{4412} into \eqref{4410}, we derive that
\begin{align*}
&(1-i)\left(j(q;q^{4})j(q;q^{2})+\alpha^{-1}j(-q;q^{4})j(-q;q^{2})\right)-
2\varphi(iq)j(\alpha;iq)\\
&=(-1+\sqrt{2}-i)j(q^2;q^4)j(-q^{6};q^{16})
-(1+\sqrt{2}+i)qj(q^2;q^4)j(-q^{14};q^{16})\\
&=\sqrt{2}(1-\alpha)j(q^2;q^4)j(-q^{6};q^{16})
-\sqrt{2}(1+\alpha)qj(q^2;q^4)j(-q^{14};q^{16})
\\&=\sqrt{2}j(q^2;q^4)j(\alpha;iq),
\end{align*}
where we use \eqref{jq} to obtain the last equality. Now we complete the proof of Entry 1.4.
\qed


\begin{thebibliography}{99}

\setlength{\itemsep}{-.8mm}
\bibitem{part1} G.E. Andrews, B.C. Berndt, Ramanujan's Lost Notebook, Part I, Springer, New York,
2005.

\bibitem{part5} G.E. Andrews, B.C. Berndt, Ramanujan's Lost Notebook, Part V, Springer, New York,
2018.

\bibitem{an} G.E. Andrews, B.C. Berndt, S.H. Chan, S. Kim, A. Malik, Four identities for third
 order mock theta functions, Nagoya Math. J., DOI: 10.1017/nmj.2018.35.

\bibitem{atkindyer}A.O.L. Atkin, P. Swinnerton-Dyer, Some properties of partitions, Proc.
London Math. Soc. (3) 4 (1954), 84--106.

\bibitem{dyson}F.J. Dyson, Some guesses in the theory of partitions, Eureka, Cambridge, 8 (1944),
10--15.

\bibitem{Gasper-Rahman-2004} G. Gasper, M. Rahman, Basic Hypergeometric Series, 2nd Ed., Cambridge Univ. Press, Cambridge, 2004.

\bibitem{hickson}D. Hickerson, A proof of the mock theta conjectures, Invent. Math. 94 (1988), 639--660.

\bibitem{hm2}D. Hickerson, E. Mortenson, Hecke-type double sums, Appell-Lerch sums, and mock theta
    functions, I, Proc. London Math. Soc. (3) 109 (2014), 382--442.

\bibitem{mo2}E. Mortenson, On three third order mock theta functions and Hecke-type double sums, Ramanujan J. 30 (2013), 279--308.

\bibitem{mo}E. Mortenson, On ranks and cranks of partitions modulo $4$ and $8$, J. Combin. Theory. Ser. A 161 (2019), 51--80.

\bibitem{ra}S. Ramanujan, The Lost Notebook and Other Unpublished Papers, Narosa Publishing House, New Delhi, 1988.

\bibitem{ye} H. Yesilyurt, Four identities related to third order mock theta functions in Ramanujan's lost
notebook, Adv. Math. 190 (2005), 278--299.

\end{thebibliography}
\end{document}